\documentclass[twocolumn]{autart_1}

\usepackage{amsmath}

\usepackage{tikz}
\usepackage{siunitx}
\usepackage[utf8]{inputenc}

\usetikzlibrary{calc}

\begin{document}
\begin{frontmatter}
\title{Non-Uniqueness of Stationary Solutions in Extremum Seeking Control}
\author[KTH]{Olle Trollberg}\ead{olletr@kth.se},
\author[KTH]{Elling W. Jacobsen}\ead{jacobsen@kth.se}

\address[KTH]{KTH Royal Institute of Technology, School of Electrical Engineering and Computer Science, Department of Automatic Control, Osquldas v. 10, 100 44 Stockholm, Sweden}

\begin{abstract}
Extremum seeking control (ESC) is a classical adaptive control method for steady-state optimization,  purely based on output feedback and without the need for a plant model.
It is well known that the extremum seeking control loop, under certain mild conditions on the controller, has a stable stationary periodic
solution in the vicinity of an extremum point of the steady-state input-output map of the plant. However, this is a local result only and 
this paper investigates whether this solution is necessarily unique given that the underlying optimization problem is convex. We first derive a necessary condition that any stationary solution 
of the ESC loop must satisfy. For plants in which the extremum point is due to a purely static nonlinearity, such as in Hammerstein or Wiener plants, the condition
involves the steady-state gradient. However, for more general plants the necessary condition involves the phase lag of the locally linearized plant, indicating
the possible existence of solutions without any relationship to optimality. Combining the derived necessary condition with 
the existence of a local solution close to the optimum, we employ the implicit function theorem and bifurcation theory to trace out branches
of stationary solutions for varying loop parameters. The focus is on solutions corresponding to limit cycles of the same period as the
forcing. We derive conditions on when these branches bifurcate, resulting in multiple stationary solutions. The results show
that cyclic fold bifurcations may exist, resulting in the existence of multiple stationary period one solutions, of which only one is related
to the optimality conditions. We illustrate the results through an example in which the conversion in a chemical reactor is optimized using ESC. 
We show that at least five stationary solutions may exist simultaneously for realistic control parameters, and that several of these solutions are stable.
One consequence of the non-uniqueness is that one in general needs to start close to the optimum to ensure convergence to the near-optimal
solution. 
\end{abstract}
\begin{keyword}
extremum seeking control; solution multiplicity; zero dynamics; bifurcations
\end{keyword}
\end{frontmatter}

\section{Introduction}
Extremum seeking control is an adaptive control method used to locate and track steady-state optima without requiring access to a process model other than an estimate of the steady-state gradient. The optimization is performed in real time and is strictly based on feedback from output measurements which makes the method applicable also in cases where a sufficiently accurate mathematical model of the process cannot be obtained. The classical perturbation based ESC method considered in this paper dates back to the early 1900's \cite{leblanc1922} and was first developed for optimization of static plants, but later applied also to dynamic plants. Most early results for dynamic plants were confined to the case of Hammerstein-Wiener like plants where linear dynamics are connected in series with a static nonlinearity, and this class of plants are still frequently considered in the literature, e.g.,  \cite{moase2012,krstic2000b}, even though it excludes many processes of interest for ESC, e.g., most chemical processes. However, in a seminal paper by Krsti\'c and Wang \cite{krstic2000} it was shown, using a combination of time-scale arguments and local analysis about the optimum, that the classical ESC method will possess a stable stationary solution in a neighborhood of the optimum also for a much wider class of dynamic plants, hence potentially widening the applicability of the method. However, existence and stability of a stationary solution is not sufficient to guarantee convergence to this solution; also uniqueness and domain of attraction need to be considered. This is exemplified in the application paper \cite{trollberg2014}, where the method is shown to possess multiple stationary solutions when applied for the optimization of a biochemical process used for wastewater treatment. The observed multiplicity implies that global convergence to the near-optimal solution cannot be guaranteed for the class of systems considered in \cite{krstic2000}. Note though that Tan, Ne{\v{s}}i{\'c}, and Mareels \cite{tan2006} have shown semi-global practical asymptotic stability of the near-optimal solution when the tuning parameters, and thereby the convergence rate, is made to approach zero. However, the observations in \cite{trollberg2014} demonstrate that results based on asymptotic time-scale separation arguments and local analysis about the optimum, such as \cite{krstic2000,tan2006}, require assumptions which can be too limiting for practical applications.

In this paper we consider the stationary solutions to the classical ESC-loop without respect to the optimality conditions and largely without resorting to asymptotic time-scale separation methods. Instead, we rely on elements of bifurcation theory to study how the stationary solutions depend on system properties and on the tuning parameters of the ESC-loop.
In Section 2, we start out by deriving necessary conditions for the stationary solutions of the loop, and show that the local phase-lag of the controlled system is critical for the existence of such solutions. In Section 3 we note that this result, when combined with the existence result of Krsti\'c and Wang \cite{krstic2000}, suggests that the optimality conditions are connected to the phase-lag of the plant. We show that the zero dynamics of the plant bifurcate at extremum points of the steady-state input-output map, and that this causes the local phase-lag to vary such that the phase-lag condition for stationarity is satisfied locally about the optimum. This hence explains the existence of a near-optimal stationary solution without resorting to time-scale separation techniques. In Section 4, we apply bifurcation theory to study the uniqueness (or lack thereof) of the near-optimal stationary solution, i.e., how solution multiplicity appears in the loop and how it is related to the tuning parameters and properties of the controlled plant. In section 5, we consider the stability of the stationary solutions, mainly to demonstrate that essentially any of the stationary solutions considered may be made stable by selection of the integral gain in the loop. Finally, the results are illustrated by means of an example in Section 6.

\section{Necessary conditions for stationarity}
In this section we derive necessary conditions for existence of periodic stationary solutions to the ESC loop. 
\subsection{The perturbation-based ESC method}
\begin{figure}[htb!]
\centering
\begin{tikzpicture}[inner sep = 0.3em, node distance =4 em]
\node[draw, circle,] at (1,1) (sum) {$\Sigma$};
\node[draw, rectangle, right of= sum,minimum width = 3em] (int) {$\dfrac{k}{s}$};
\node[draw, rectangle, right of= int,node distance=5em] (FL) {$F_L(s)$};
\node[draw, circle, inner sep=0.3em,right of = FL]  (prod) {$\Pi$};
\node[draw, rectangle, above of= prod,node distance = 4em] (FH) {$F_H(s)$};
\node[draw, rectangle, node distance = 5em] at ($($(int)!0.5!(FL)$)+(0,6em)$)(sys) {$\begin{aligned}\dot x &= f(x,u) \\ y &= h(x)\end{aligned}$};
\draw[-stealth] (sum) ++ (0,-1) -- node[pos=0.1,right]{$a\sin(\omega t)$} (sum); 
\draw[-stealth] (prod) ++ (0,-1) -- node[pos=0.1,left]{$\sin(\omega t)$} (prod); 
\draw[-stealth] (prod) -- (FL);
\draw[-stealth] (FL) --node[above] {$\xi$} (int);
\draw[-stealth] (int) -- node[above] {$\hat u$}(sum);
\draw[-stealth] (sum) |- node[above,pos = 0.75] {$u$}(sys);
\draw[-stealth] (sys) -|node[pos=0.25, above]{$y$} (FH);
\draw[-stealth] (FH) -- node[left]{$y-\eta$} (prod);
\end{tikzpicture}
\caption{The classic perturbation based ESC loop}
\label{fig:esc_loop}
\end{figure}
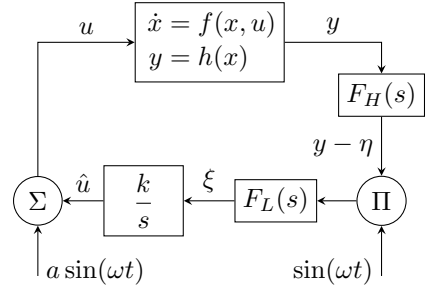

Figure~\ref{fig:esc_loop} illustrates the classical ESC loop and defines the various signals used in the scheme. The basic operation may tentatively be described as follows: A sinusoidal perturbation is added to the nominal input $\hat u$ in order to excite the system to reveal gradient information. For static systems, the amplitude of the response to the perturbation in the output $y$ becomes proportional to the local gradient. The amplitude is then extracted by means of high-pass filtering through $F_H(s)$ and demodulation by multiplication with a signal of the same frequency as the perturbation. The demodulation introduces high-frequency byproducts which are attenuated by a low-pass filter $F_L(s)$. The output $\xi$ of the low-pass filter is then approximately proportional to the local gradient, and by closing the loop with an integrator with the correct sign, the loop is driven towards a point of zero gradient, i.e., the optimum. For dynamic systems, this tentative argument breaks down, but Krsti\'c and Wang \cite{krstic2000} have nevertheless rigorously shown that a stable near-optimal periodic stationary solution exists under certain conditions also for dynamic plants.

In this paper we follow Krsti\'c and Wang \cite{krstic2000} and consider plants of the form
\begin{equation}
\begin{aligned}
\dot x =& f(x,u), && x\in \mathbf R^n, u \in \mathbf R \\
y =& h(x), &&  y \in \mathbf R
\end{aligned}
\label{eq:sys}
\end{equation}
where we assume that \eqref{eq:sys} is open-loop stable and that the steady-states are parametrized by the input through a function $l:\mathbf R \rightarrow \mathbf R^n$ such that
\begin{equation}
f(x,u) = 0 \text{ if and only if } x = l(u).
\end{equation}
For lack of better terminology, this class of plants is referred to as plants with `general' dynamics.
The objective of the ESC loop is to maximize\footnote{We consider maximization in this paper, but the results trivially also hold for minimization.} the composite function
\begin{equation}
J(u) = h \circ l (u)
\end{equation}
which corresponds to the steady-state input-output map, sometimes referred to as the equilibrium map, of the plant.
All functions $f$, $h$, $l$, and $J$ are assumed to be sufficiently smooth for all necessary derivatives to exist.
If first order filters are used, $F_H(s) = s/(s+\omega_h)$ and $F_L(s) = \omega_l/(s+\omega_l)$, then the complete dynamics of the closed loop may be expressed as
\begin{equation}
\begin{aligned}
\dot x &= f(x,\hat u + a\sin(\omega t)) \\
y &= h(x)\\
\dot \eta &= \omega_h(y - \eta)\\
\dot \xi &= \omega_l((y-\eta)\sin(\omega t) - \xi)\\
\dot{\hat u} &= k \xi,
\end{aligned}
\label{eq:cl}
\end{equation}
but the results derived in this paper does not depend critically on the particular choice of filters.

\subsection{Stationary solutions of the ESC loop}
Since the closed loop includes a periodic forcing, see Figure~\ref{fig:esc_loop} and \eqref{eq:cl}, it is unlikely that any steady-state solution exists. In this paper, we will thus consider periodic stationary solutions with the same fundamental period as the forcing, i.e., $T=2\pi/\omega$. Other types of stationary solutions, such as double period and chaotic solutions, may in principle also occur, but then an underlying period-one solution typically still exists, although it may be unstable. 

If we let $z = [x\ \ \hat u \ \ \xi \ \ \eta]^T$ be an extended state vector with both the system and controller states, we may represent the closed loop system  \eqref{eq:cl} as
\begin{equation}
\dot z = F(z,t).
\end{equation}

Assume now that the system is operated at a periodic stationary solution with period $T=2\pi/\omega$.
Due to the periodicity and the imposed smoothness, it must hold that 
\begin{equation}
\int_t^{t+T}F(z(\tau),\tau)d\tau = 0.
\label{eq:periodicity}
\end{equation}
We may also expand the states $z_i$ and their derivatives $\dot z_i$ by Fourier series. 
Let the Fourier series of $\hat u(t)$ be
\begin{equation}
\hat u(t) = \sum_{k=-\infty}^\infty c^{\hat u}_ke^{ik\omega t}.
\end{equation}
where $c_k^ \cdot$ are complex coefficients.
The input to the plant, $u(t)$, may then be expressed as
\begin{equation}
u(t) = \sum_{k=-\infty}^\infty c^{\hat u}_ke^{ik\omega t} + \frac{a}{2i}(e^{i\omega t} - e^{-i\omega t})
\end{equation}
where the additional terms are introduced by the external perturbation.
Assume now that the amplitude of the input is small such that the plant may be locally approximated as linear, and let $\bar u = c^{\hat u}_0$ be the average value of the input over one period. Then $u = \bar u$, $x=l(\bar u)$ may be considered the operating point associated with the stationary solution. 
Let
\begin{equation}
G_{\bar u}(s) = \frac{\partial h}{\partial x}\left (sI - \frac{\partial f}{\partial x} \right)^{-1}\frac{\partial f}{\partial u}
\label{eq:linearization}
\end{equation}
be the transfer function of the locally linearized plant, i.e., the Jacobians are evaluated at $u =\bar u, x = l(\bar u)$. The plant-response to low-amplitude inputs may then be described as
\small
\begin{align}
y(t) &=  \sum_{k=-\infty}^\infty c^{y}_ke^{ik\omega t} \\ 
 = &\sum_{k=-\infty}^\infty G_{\bar u}(ik\omega)c^{\hat u}_ke^{ik\omega t} + \frac{a(G_{\bar u}(i\omega)e^{i\omega t} -G_{\bar u}(-i\omega) e^{-i\omega t})}{2i}.
\end{align}
\normalsize
For the high-pass filter we have
\begin{equation}
\dot \eta(t) = \omega_h(y(t)-\eta(t)) = \omega_h\sum_{k=-\infty}^\infty c^{\dot\eta}_ke^{ik\omega t}
\label{eq:dotetafourier}
\end{equation}
where 
\begin{equation}
c_k^{\dot\eta} = F_H(ik\omega)c_k^y.
\label{eq:highpass_coeff}
\end{equation}
Furthermore, since \eqref{eq:periodicity} must hold for each component in $z$, we have
\begin{equation}
\int_t^{t+T}\dot{\hat u} d\tau = k\int_t^{t+T}\xi d\tau = 0
\label{eq:dk}
\end{equation}
and
\begin{equation}
\int_t^{t+T}\dot\xi d\tau = \omega_l\int_t^{t+T}((y-\eta)\sin(\omega \tau) - \xi)d\tau = 0
\label{eq:dxi}
\end{equation}
for any stationary solution with period $T$.

By substituting \eqref{eq:dk} and \eqref{eq:dotetafourier} into \eqref{eq:dxi}, and using Euler's identity on the sine term, we obtain
\begin{equation}
\sum_{k=-\infty}^\infty\int_t^{t+T} \frac{c^{\dot\eta}_k(e^{i(k+1)\omega t}-e^{i(k-1)\omega t})}{2i} d\tau = 0.
\end{equation}
Most terms integrate to zero, leaving only
\begin{equation}
\int_t^{t+T} \frac{c^{\dot\eta}_{-1}-c^{\dot\eta}_{1}}{2i} d\tau = 0
\end{equation}
which implies 
\begin{equation}
c^{\dot\eta}_{-1} = c^{\dot\eta}_{1}.
\end{equation}
Since $\dot \eta$ is real, we have symmetry in the Fourier series coefficients such that $c_{-1}^{\dot \eta}$ is the complex conjugate of $c_{1}^{\dot \eta}$.
We thus conclude that we must have
\begin{equation}
\operatorname{Im}\{c^{\dot\eta}_{1}\} = 0
\end{equation}
at any low-amplitude stationary solution of period $T$.
By substituting \eqref{eq:highpass_coeff} into this condition, we see that we must have
\begin{equation}
\operatorname{Im}\left\{F_H(i\omega)G_{\bar u}(i\omega)\left(c_1^{\hat u}+\frac{a}{2i}\right)\right\} = 0
\label{eq:condition}
\end{equation}
for any periodic stationary solution with period $T=2\pi/\omega$ and low amplitude. 

\subsection{Interpretation as a phase condition}
It is useful to consider \eqref{eq:condition} as a condition on the local properties of the controlled system. For this purpose, assume that the break-off frequency of the low-pass filter is much smaller than $\omega$, i.e., $\omega_l \ll \omega$ such that $c_1^{\hat u} \ll a$. Note that we normally expect $c_1^{\hat u} \ll a$ to hold even without this assumption for conservatively tuned ESC loops since harmonics in the loop are damped both by the low-pass filter and the integral controller. Nevertheless, using this assumption we may neglect the $O(c_1)$ term and simplify the stationarity condition to
\begin{equation}
\operatorname{Im}\left\{F_H(i\omega)G_{\bar u}(i\omega)\frac{a}{2i} \right \} = 0
\end{equation}
which implies that
\begin{equation}
\operatorname{Re}\left\{F_H(i\omega)G_{\bar u}(i\omega)\right\} = 0.
\label{eq:stat_cond_lp}
\end{equation}
This condition is trivially satisfied when $G_{\bar u}(i\omega)=0$. In cases where $G_{\bar u} (i\omega)\neq0$, condition \eqref{eq:stat_cond_lp} may instead be interpreted as a condition on the local phase-lag of the controlled plant, i.e.,
\begin{equation}
\angle G_{\bar u} (i\omega) = \frac{\pi}{2} - \angle F_H(i\omega) + n\pi, \quad n\in \mathbf Z
\label{eq:phase_cond}
\end{equation}
where $\mathbf Z$ is the set of integers and $\angle F_H(i\omega)$ is a constant determined by the tuning parameters used.

\section{Connecting optimality and stationarity}

Above, it was shown that the periodic stationary solutions of the loop can be characterized by condition \eqref{eq:stat_cond_lp}. For systems where the steady-state optimum corresponds to an extremum in a static nonlinearity in series with dynamic elements, such as in the case of Hammerstein or Wiener systems, the plant will be locally invariant at the optimum, i.e., $G_{\bar{u}}(i\omega)\equiv 0$, and condition \eqref{eq:stat_cond_lp} will hence be trivially satisfied there. However, in the general case, plants will have a dynamic response also at the steady-state optimum and hence $G_{\bar{u}}(i\omega)\not \equiv 0$ at any steady-state operating point. Nevertheless, given that Krsti\`c and Wang \cite{krstic2000} have proven existence of a near-optimal stationary solution in the general case, condition \eqref{eq:stat_cond_lp} must hold at some point in a neighborhood of the optimum, although it may not be immediately clear why this should be the case. Since the condition is dynamic in nature, previous results on existence, e.g., \cite{krstic2000,tan2006}, provide limited insight as they rely on the use of asymptotic arguments to essentially reduce the controlled plant to a static map such that all information on the local dynamic properties is lost. In the following, we will approach the problem in a fully dynamic setting.

Bifurcation theory provides a link between the stability of a dynamical system and the branching behavior of its
stationary solutions; solution branches meet where eigenvalues of the linearized dynamics cross the stability boundary
\cite{guckenheimer2013}. For the case of static bifurcations, it implies that certain dynamic properties 
can be predicted from steady-state information about the system only, e.g., a singularity in the steady-state input-output map implies that an eigenvalue crosses the imaginary axis at that point and at least one of the steady-state
branches emerging from the singularity will be unstable.  Similarily, it is natural to consider the stability of the inverse dynamics (zero-dynamics) with respect to the properties of the steady-state input-output map. Clearly, an extremum point in the steady-state input-output map corresponds to a singularity in the inverse map, i.e., the steady-state output-input map. This indicates that extremum points in the static map are connected to bifurcations in the plant zero-dynamics. Bifurcations in the zero-dynamics in turn implies that a zero crosses the stability boundary with large local variations in the phase-lag as a consequence. This ensures that the phase condition \eqref{eq:phase_cond} is satisfied locally about extrema in the steady-state input-output map. Below, we formalize this argument.

\subsection{Bifurcations of the zero dynamics}
Here we consider single-input single-output nonlinear dynamic systems described by a set of ordinary differential and algebraic equations on input-affine form
\begin{equation}
\begin{aligned}
\dot{x}=&f(x)+g(x)u, \ x \in \mathbf{R}^n, \ u \in \mathbf{R}\\
y=&h(x), \ y \in \mathbf{R} \label{zero:eq:mod}
\end{aligned}
\end{equation}

Note that this only covers a subset of the systems described in Section 2.1, but that the main results derived below apply also to systems that can not be written on input-affine form. However, in such cases, it may be challenging to find explicit expressions for the zero-dynamics whereby the derivations become more involved and are therefore not included here. 

The zero dynamics of a system correspond to the state dynamics when the output $y$ is forced to be zero or, more generally, constant \cite{sastry2013nonlinear}. To determine the zero dynamics of the system \eqref{zero:eq:mod}, we introduce a state transformation\footnote{This is always possible for systems on the form \eqref{zero:eq:mod} under mild assumptions \cite{sastry2013nonlinear}, and this is the main reason for considering input-affine systems. However, note that the zero dynamics usually are well defined also in cases where it is difficult or not possible to transform the problem into normal form such that the zero dynamics become explicit.} $z=\phi(x)$ to obtain the normal form
\begin{equation}
\begin{aligned}
\dot{z}_i &= z_{i+1}, \quad i=1,\ldots,r-1 \\
\dot{z}_r &= b(\chi,\psi)+a(\chi,\psi)u \\
\dot{\psi} &= q(\chi,\psi)\\
y &= z_1 
\end{aligned}
\end{equation}
where $\chi=\left[ z_1 \ z_2  \ldots  z_r \right]$, $\psi=\left[z_{r+1} \ \ldots z_n\right]$, and $r$ is the relative degree of the system. The zero dynamics are then given by the dynamics of the $n-r$ dimensional state $\psi$ when the $r$-dimensional state $\chi$ is forced to be zero by means of the control input $u$, i.e., 
\begin{equation} 
\dot{\psi}=q(0,\psi). \label{zero:eq:zd}
\end{equation} 

\subsection{Static bifurcations of the zero dynamics}
We are here concerned with consequences of static bifurcations of the zero dynamics \eqref{zero:eq:zd} as we move the operating point along the equilibrium manifold of (\ref{zero:eq:mod}). A bifurcation of the zero dynamics occurs when eigenvalues of $q(0,\psi)$, linearized about a point on the equilibrium map, cross the imaginary axis. As shown in \cite{Isidori89}, the linear approximation of the zero dynamics at an equilibrium point equals the zero dynamics of the linearized system at the same equilibrium. That is, eigenvalues of the linearized zero dynamics coincide with the zeros of the linearized dynamics of the open-loop system \eqref{zero:eq:mod}. Bifurcations can hence be determined from consideration of the transmission zeros of
\begin{equation}
\begin{aligned}
\dot{x}&=Ax(t)+Bu(t) \\
y(t)&=Cx(t) \label{zero:eq:lin}
\end{aligned}
\end{equation}
where $(A,B,C)$ is the linear approximation of \eqref{zero:eq:mod} around a given steady state. 

The transmission zeros of the linearized system \eqref{zero:eq:lin} can be determined from the rank of the Rosenbrock system-matrix
\begin{equation}
M=\left( \begin{matrix} zI-A\, \ & -B \\ C & 0\end{matrix} \right). \label{zero:rank1}
\end{equation}
The transmission zeros are the values of $z$ such that the rank of $M$ is less than the normal rank \cite{rosenbrock1970state}. Using Schur's determinant formula we get
\begin{equation}
\det(M)=\det (zI-A)\det C(zI-A)^{-1}B=0 \label{eqz}
\end{equation}
and we can hence conclude that $z$ is a zero if $\det C(zI-A)^{-1}B=0$ and $z$ is not an eigenvalue of $A$. The latter condition rules out pole-zero cancellations. At a bifurcation point of the zero dynamics, at least one zero will have real part $\Re\{z\}=0$.

A static bifurcation of the zero-dynamics, i.e., generally a fold bifurcation, implies that $z=0$. From (\ref{eqz}), this condition translates into $CA^{-1}B=0$ which as expected corresponds to a zero steady-state gain $G(0)=0$ from input to output. However, to be a bifurcation point, a transversality condition also needs to be satisfied, i.e., the zero must cross the imaginary axis as the equilibrium point is varied. For this purpose, consider the
MacLaurin series of $G(s)=C(sI-A)^{-1}B$ 
\begin{equation}
G(s)=\Sigma_{i=0}^\infty c_is^i
\end{equation}
where $c_i=CA^{-1-i}B$. For small non-zero $s$, we can neglect higher order terms which implies that the zero close to $s=0$ is given by
\begin{equation}
z=-\frac{c_0}{c_1}=-\frac{CA^{-1}B}{CA^{-2}B}.
\end{equation}
If we assume that the zero at $z=0$ has multiplicity one, i.e., only a single zero moves through the origin,  then $CA^{-2}B$ must be non-zero, and we find that $CA^{-1}B=G(0)$ changes sign as the zero changes sign. Thus, a static bifurcation of the zero dynamics, corresponding to a real zero crossing the imaginary axis, implies that the local steady-state gain changes sign. This then corresponds to an extremum point in the steady-state input-output map. 

Our primary concern here is whether the converse of the above result is true, i.e., whether an extremum point in the equilibrium map implies a static bifurcation of the zero dynamics. At an extremum point we have $CA^{-1}B=0$ and we note from the MacLaurin series above that $z=0$ is then a transmission zero of $G(s)$ unless also $CA^{-i}B$, $\forall i>1$ are also all identically zero. The latter case corresponds to having $G(s) \equiv 0$ at the extremum point, and this is indeed possible if the zero gain is due to a static nonlinearity, as in Wiener and Hammerstein models. However, when the nonlinearity causing the extremum point is not static but inherent in the state dynamics, then the system will display a transient response also when operated at the extremum point, i.e.,
$G(s) \not \equiv 0$. Then $G(0)=0$ implies that a zero exist at $z=0$ and the change in the sign of $G(0)$ at the extremum point implies that the transversality condition will be satisfied. Hence, an extremum point in the equilibrium map will correspond to a static bifurcation of the zero dynamics for systems satisfying $G(s) \not \equiv 0$ at the extremum point.

We remark that the above results do not imply that at least one solution has unstable zero dynamics in the case of input multiplicity, as is sometimes claimed e.g., \cite{sistu95}. The main reason for this is that transmission zeros may move between the complex LHP and RHP either through the imaginary axis or through infinity, and the latter case does not correspond to a bifurcation and has no effect on the sign of  the  steady-state gain. Thus, all we can conclude is that a static bifurcation of the zero dynamics implies an extremum point in the equilibrium map.

\subsection{Satisfaction of phase condition due to crossing zero}

From the above we conclude that a real transmission zero is crossing the imaginary axis at an extremum point in the steady-state input-output map (unless the steady-state optimum is due to a purely static relationship). We here consider the implications of this crossing for the process dynamics and how this relates to the stationarity condition \eqref{eq:phase_cond}.

\begin{figure}[h!]
\centering
\includegraphics[width=\columnwidth]{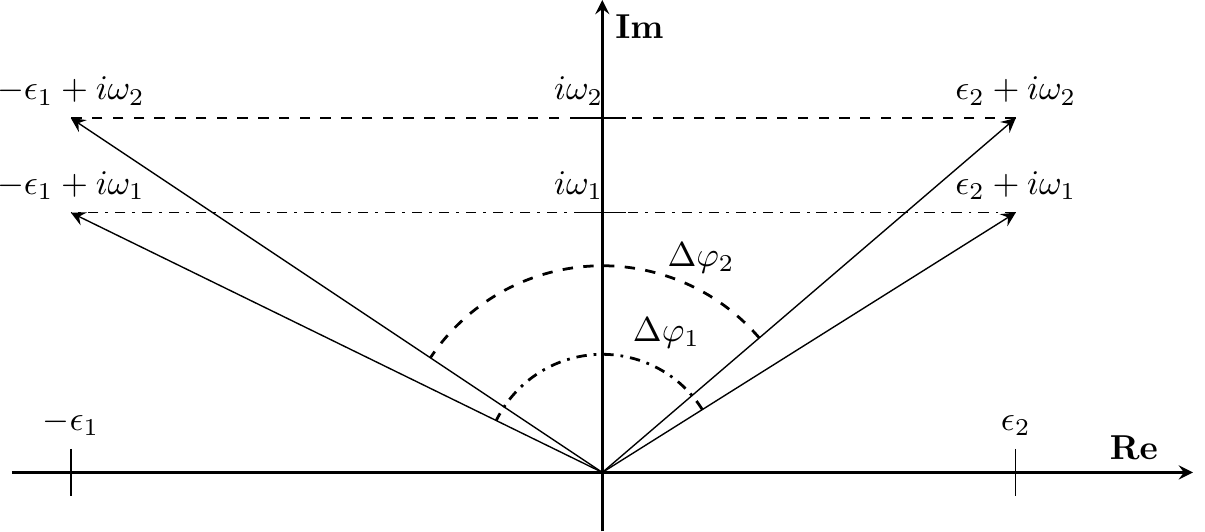}
\caption{Change in phase contribution of varying real transmission zero $z$ in a linear system on the form $G(s)=(s+z)G_0(s)$. $\Delta \varphi_i$ denotes the change in the phase-lag
of the system for the fixed frequencies $\omega_i$ as $z$ varies over the interval $(-\epsilon_1,\epsilon_2)$, $\epsilon_1,\epsilon_2>0$.}
\label{zero:fig:shift}
\end{figure}

The fact that the transfer-function $G(s)$ at an extremum point contains a zero at the origin implies partly that the steady-state gain is zero at such points, and partly that the phase-lag is $\pm \pi/2$ rad at $\omega=0$. That the zero crosses the imaginary axis as the operating point, i.e., the point of linearization, moves past the extremum point, implies that the steady-state gain remains small close to the extremum point while the phase-lag at $\omega=0$ change by $\SI{\pi}{\radian}$, centered about the phase-lag $\pm\SI{\pi/2}{\radian}$. For low nonzero frequencies $\omega$, the phase-characteristics of $G(s)$ near an extremum point remains, see Figure~\ref{zero:fig:shift}. At an extremum point, the phase-lag contributed by the zero to the local frequency response function $G(i\omega)$ will be $\pm\SI{\pi/2}{\radian}$. As the operating point is varied locally past the optimum, the zero will vary over some interval $z\in[-\epsilon_1,\epsilon_2]$ with $\epsilon_1,\epsilon_2>0$. The phase contribution to $G(i\omega)$ from the zero hence varies with $\Delta\varphi$, approximately centered about $\pm\SI{\pi/2}{\radian}$, see Figure~\ref{zero:fig:shift}. For sufficiently low frequencies we have $\Delta\varphi\approx\SI{\pi}{\radian}$. This phase-variation guarantees satisfaction of the phase condition \eqref{eq:phase_cond} locally about an extremum point, at least for sufficiently low frequencies.

To see how the variations in the phase-lag and the crossing zero relates to the satisfaction of \eqref{eq:phase_cond}, we first make the assumption that the high-pass filter $F_H(i\omega)$ is defined with a break-off frequency $\omega_h = \alpha \omega$ with $\alpha < 1$ such that the phase lag due to this filter is independent of $\omega$. Note that we require $\omega_h < \omega$ for the filter to operate as intended, so this assumption is quite natural.
Next, we write the plant as $G_{\bar u}(s)=(s+z_{\bar u})G_{\bar u}^0(s)$ where $z_{\bar u}$ is the zero crossing the imaginary axis, and $G_{\bar u}^0(s)$ collects the remaining poles and zeros of $G(s)$. We let the input ${\bar u}$ represent the operating point $(x,{\bar u})=(l({\bar u}),{\bar u})$ about which we linearize and use subscript ${\bar u}$ to indicate the dependence on the operating point. If we for simplicity neglect the multiple of $\pi$ and only consider the positive case\footnote{Neglecting these factors only affects the sign in the following expressions, and we leave it to the reader to fill in the remaining cases.}, then the phase condition \eqref{eq:phase_cond} may be expressed as
\begin{equation}
\angle(i\omega + z_{\bar u}) + \angle G_{\bar u}^0(i\omega) + \angle F_H(i\omega) = \frac{\pi}{2}
\end{equation}
or equivalently as
\begin{equation}
\tan^{-1} \frac{\omega}{z_{\bar u}} = \frac{\pi}{2} - \angle G_{\bar u}^0(i\omega) - \angle F_H(i\omega).
\end{equation}
If we solve for $z_{\bar u}$ we get
\begin{equation}
z_{\bar u} = \frac{\omega}{\tan \left(\frac{\pi}{2} - \angle G_{\bar u}^0(i\omega) - \angle F_H(i\omega)\right)}. \label{eq:zu}
\end{equation}
Essentially, we have above translated the phase condition \eqref{eq:phase_cond} into a condition on the crossing zero. Since the zero will move continuously through some interval containing the origin as the operating point is varied in a neighborhood of an extremum point, this condition will be satisfied near such points given that the right hand side of \eqref{eq:zu} stays sufficiently close to zero. For low enough frequencies $\omega$, this will be the case as we show below.

Since $G_{\bar u}^0(i\omega)$ will not contain any poles or zeros at the origin, we either have $\angle G_{\bar u}^0(i\omega) \rightarrow 0$ or $\angle G_{\bar u}^0(i\omega) \rightarrow \pi$ as $\omega \rightarrow 0$. Together with the assumption on $F_H(s)$, this implies that the denominator in \eqref{eq:zu} will approach a nonzero constant when $\omega\rightarrow0$. It hence follows that the right hand side of \eqref{eq:zu} approach zero as $\omega \rightarrow 0$. This implies that condition \eqref{eq:zu} will be satisfied for sufficiently low frequencies. 
To conclude, we should expect a stationary solution to exist near an extremum in the steady-state input-output map.

We next make use of \eqref{eq:zu} to consider the deviation of a stationary solution from the extremum point. Since $z_{\bar u}$ will vary continuously with the operating point, we can locally about an extremum point ${\bar u}^*$ use the linear approximation
\begin{equation}
z_{\bar u} \approx \underbrace{z_{{\bar u}^*}}_{=0} + ({\bar u}-{\bar u}^*)\frac{dz_{\bar u}}{d{\bar u}}
\end{equation}
from which we may conclude that the deviation from the optimum may be approximated by
\begin{equation}
{\bar u}-{\bar u}^*\approx \frac{\omega}{\tan \left(\frac{\pi}{2} - \angle G_{\bar u}^0(i\omega) - \angle F_H(i\omega)\right) \frac{dz_{\bar u}}{d{\bar u}}}.
\end{equation}

It is interesting to note that the above results are based purely on the phase-lag properties of the plant and filters, and that a steady-state optimum hence in principle may be located using phase-information only. This opens up for a novel field of ESC-algorithms based on phase-estimation and control. However, we do not pursue this topic further here other than noting that phase-locked loops \cite{gardner2005phaselock} may be of central interest in such an approach.
\section{Solution multiplicity in ESC}
In the previous two sections we established a necessary condition for stationarity and showed that this condition will be satisfied in a neighborhood of an extremum point in the ESC objective function. However, the condition may also be satisfied at operating points unrelated to the steady-state optimum, hence suggesting that the ESC-loop may suffer from existence of multiple stationary solutions of which some may be unrelated to the optimality conditions. In fact, such multiplicity has previously been reported in the literature \cite{trollberg2014}. In this section, we employ bifurcation theory to consider how such solution multiplicity may appear in the ESC loop and how it is related to properties of the controlled system.

Krsti\'c and Wang \cite{krstic2000} have established that the tuning parameters of the ESC-loop may be chosen such that a stable stationary periodic solution exists in a neighborhood of the steady-state optimum. Whenever a solution exists, the implicit function theorem implies that a branch of qualitatively similar solutions also exists when the loop parameters are varied (at least locally). This branch may generally be continued until a possible bifurcation point is reached. At a bifurcation point, at least two solution branches meet and the stability of the solutions change. Existence of bifurcation points are hence of prime interest for us since they imply local solution multiplicity. Since we are mainly interested in existence of multiple period-one solutions, we here focus on existence of cyclic fold bifurcations where two such branches meet \cite{guckenheimer2013}.

Since an analytic treatment for the full loop described by \eqref{eq:cl} is hard in a general setting, 
we will approach the problem via a proxy, namely the necessary condition \eqref{eq:stat_cond_lp} derived above. Being necessary, this condition must be satisfied at all low-amplitude periodic solutions of the loop with a period equal to the applied forcing.
Hence, whenever there exist such a branch of periodic solutions, we can look at the necessary condition locally, and if this goes through a fold bifurcation, so must the periodic branch (more accurately, the periodic branch will go through a cyclic fold bifurcation) since it cannot exist past the fold in the necessary condition. Note that the periodic solution-branch may exhibit other types of bifurcations not detected in the necessary condition, but the original period-one branch will in that case persist, possibly alongside other solution-branches, as long as the solution does not pass through any singularity like a fold bifurcation. A bifurcation of the necessary condition is hence sufficient for the existence of multiple stationary solutions of ESC.

In order to make the analysis below tractable, we make a few assumptions on the loop tuning parameters. We will assume that the break-off frequency of the low-pass filter satisfies $\omega_l \ll \omega$ such that $c_1^{\hat u} \ll a$ and can be neglected, i.e., we assume that the harmonics are damped out completely such that we can use the simplified necessary condition \eqref{eq:stat_cond_lp}. Such a low break-off frequency would cause the convergence rate to become impractically low, but it significantly simplifies the analysis by enabling us to focus on system properties and neglect ``standing waves'' in the loop. We will also assume that the perturbation amplitude $a$ is small such that the system locally may be approximated as linear which was also assumed in the derivation of \eqref{eq:stat_cond_lp}. Finally, we restrict our analysis to the parameter $\omega$, i.e., the perturbation frequency. This choice is natural since the conditions we consider are strongly related to the local frequency-response of the system. However, note that in the generic case, a bifurcation point for one parameter will also be a bifurcation point when other parameters are considered.

\subsection{Conditions for existence of a fold bifurcation}

To make clear the dependence on the operating point and the bifurcation parameter $\omega$, let condition
 \eqref{eq:stat_cond_lp} be written
\begin{equation}
C(\bar u,\omega) = \operatorname{Re}\left\{F_H(i\omega)G_{\bar u}(i\omega)\right\} = 0
\label{eq:param}
\end{equation}
where we let $\bar u$ represent the operating point (average input over one period), and where the local transfer function $G_{\bar u}(s)$ is defined by \eqref{eq:linearization}. Conditions for a fold bifurcation (or turning point) in an algebraic relation such as \eqref{eq:param} are \cite{guckenheimer2013,seydel2010}
\begin{subequations}
\begin{align}
C(\bar u,\omega) = 0, \label{eq:conda}\\
\frac{\partial C}{\partial \bar u} = 0, \label{eq:condb} \\
\frac{\partial C}{\partial \omega}\neq 0,\label{eq:condc}\\
\frac{\partial^2 C}{\partial \bar u^2} \neq 0,\label{eq:condd}
\end{align}
\end{subequations}
Here \eqref{eq:conda}, \eqref{eq:condb}, and \eqref{eq:condc} ensure that the solution branch is locally perpendicular to the parameter-axis, and \eqref{eq:condd} ensures that the branch turns back (cf. \cite[Definition 2.8, p.74]{seydel2010}). The last condition  \eqref{eq:condd} is frequently replaced by other conditions aimed at avoiding degenerate situations where the candidate point is an inflection point at which the solution branch is perpendicular to the parameter-axis but does not turn back. Such points are \emph{structurally unstable} and will for small perturbations of the problem decompose into two fold bifurcations or no bifurcation at all \cite{seydel2010}. Since \eqref{eq:condd} is mainly there to avoid these pathological cases, it is of less practical interest and we will focus our discussion on the three first conditions.

The first condition \eqref{eq:conda} simply states that the bifurcation point is part of the solution branch.
According to the second condition \eqref{eq:condb}, the partial derivative of $C$ with respect to $\bar u$ has to be zero. If we expand this derivative we get
\begin{equation}
\frac{\partial C}{\partial \bar u} = \operatorname{Re}\left\{F_H(i\omega)\frac{\partial G_{\bar u}(i\omega)}{\partial \bar u}\right\} = 0.
\end{equation}
This is true if either 
\begin{equation}
\frac{\partial G_{\bar u}(i\omega)}{\partial \bar u}=0,
\end{equation}
or if 
\begin{equation}
\angle F_H(i\omega) + \angle \frac{\partial G_{\bar u}(i\omega)}{\partial \bar u} = \frac{\pi}{2} + n\pi
\end{equation}
for $n\in\mathbf Z$.

The third condition \eqref{eq:condc} may be thought of as a continuity/smoothness condition. At the bifurcation point, the implicit function theorem fails and $\bar u$ cannot locally be considered a function of $\omega$. However, \eqref{eq:condc} implies that $\omega$ instead may be considered a function of $\bar u$. In other words, the two solution branches that meet at the fold bifurcation together form a smooth continuous solution curve in the $(\bar u, \omega)$ plane, and this curve may be uniquely continued through the bifurcation. This effectively excludes other codimension-one bifurcations such as pitchfork or transcritical bifurcations.
If we expand \eqref{eq:condc} and simplify, we find that it is equivalent to
\begin{equation}
\operatorname{Im}\left\{\frac{\partial F_H(i\omega)}{\partial s}G_{\bar u}(i\omega) + F_H(i\omega)\frac{G_{\bar u}(i\omega)}{\partial s} \right\} \neq 0
\end{equation}
where $s$ is a complex number.
When $G_{\bar u}(s)\not \equiv 0$, this is generally satisfied except possibly at singular points. That such a singularity would coincide with the other conditions for a bifurcation is unlikely and we will not pursue this pathological case further here.

\subsection{Geometric interpretation of the bifurcation conditions}
\begin{figure}[htb!]
	\centering
	\includegraphics[]{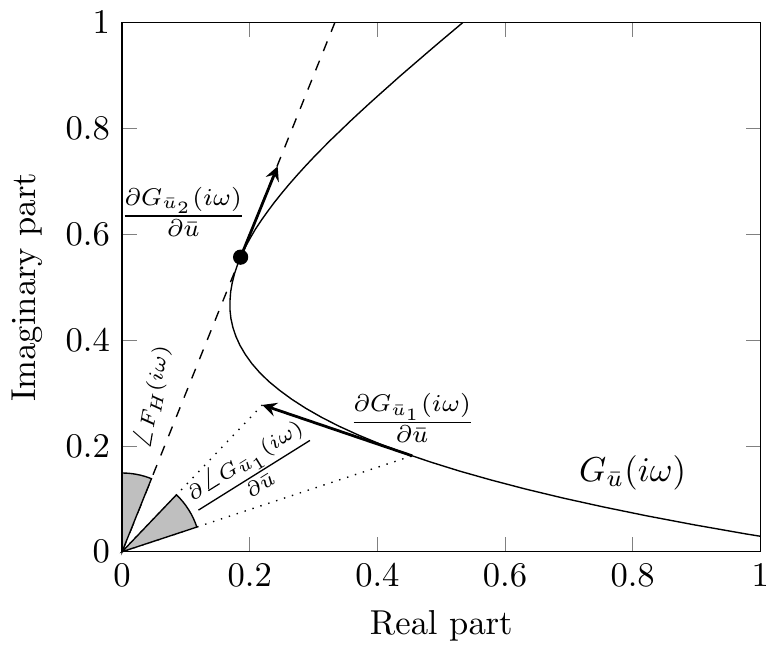}
	\caption{Solid line corresponds to $G_{\bar u}(i\omega)$ as $\bar u$ is varied and $\omega$ is fixed. Dashed ray corresponds to the argument of $F_H(i\omega)$. }
	\label{fig:geometric_interpretation}
\end{figure}

The conditions for a fold bifurcation discussed above may be interpreted geometrically, see Figure~\ref{fig:geometric_interpretation}. The solid curve represents the frequency response as a complex number in the imaginary plane for a fixed perturbation frequency $\omega$ as the operating point $\bar u$ is varied. Each point on the curve thus corresponds to a given $\bar u$. The dashed ray represents the argument of $F_H(i\omega)$. Any complex number on the ray will hence be taken to the imaginary axis when multiplied with $F_H(i\omega)$. In particular, the necessary condition for stationarity \eqref{eq:conda} is satisfied at points where the ray intersects the solid curve since the product $F_H(i\omega)G_{\bar u}(i\omega)$ then will be purely imaginary.
At any given point on the solid curve, we may evaluate $\partial G_{\bar u}/\partial \bar u$, and this then corresponds to a tangent vector of the curve. The bifurcation conditions \eqref{eq:conda}, \eqref{eq:condb} and \eqref{eq:condd} may then be interpreted as a point where the curve touches the ray tangentially without crossing. When \eqref{eq:condc} is satisfied, small perturbations of the bifurcation parameter $\omega$ will cause the curve and ray to move such that the curve either does not intersect the ray at all, or intersects the ray at two separate points, thus giving rise to solution multiplicity.

Consider next the argument of the locally linearized system $\angle G_{\bar u}(i\omega)$. 
The relation between the argument ${\partial \angle G_{\bar u}(i\omega)}/{\partial \bar u}$ and the tangent-vector ${\partial G_{\bar u}(i\omega)}/{\partial \bar u}$ is also illustrated in Figure~\ref{fig:geometric_interpretation}. Clearly, at the bifurcation point, it must hold that
\begin{equation}
\frac{\partial \angle G_{\bar u}(i\omega)}{\partial \bar u}  = 0.
\end{equation}
Hence, if we interpret the stationarity condition as a condition on the local phase-lag, then a fold bifurcation corresponds to a point where the phase-condition \eqref{eq:phase_cond} is satisfied and where the phase-lag as a function of the operating point $\bar u$ is at an extreme point. We may hence expect fold bifurcations to appear in systems where the phase-lag at the perturbation frequency varies in a non-monotone manner with the operating point.

Before moving on, we remark that the results derived in this section rely on the assumptions that $\omega_l \ll \omega$ and that $a$ is sufficiently small. Nevertheless, we believe that the results are applicable also in a more general setting, although this is difficult to prove analytically. This is supported by the results of our numerical example at the end of the paper where the amplitude $a$ and break-off frequency $\omega_l$ are non-negligible.


\section{Stability of stationary solutions}

In this section we briefly discuss how stability of the stationary solutions may be determined in order to demonstrate that essentially any of the stationary solutions discussed above may be made stable by selection of an appropriate integral gain $k$ in \eqref{eq:cl}. A rigorous stability analysis is out of scope since our main focus here is existence of multiple stationary solutions as such. 

Consider the closed-loop system \eqref{eq:cl} and assume that $k$ is small compared to the other parameters. Then we may apply singular perturbations to separate the integral state $\hat u$ into a reduced model, and the rest of the dynamics into a boundary layer model (cf. \cite{khalil2000}). The boundary layer model is stable by assumption\footnote{Note however that the origin is not exponentially stable; the perturbation will generally cause the boundary layer model to enter a periodic orbit, and this may prevent the direct application of standard results such as Tikhonov's theorem.} (series connection of stable systems, and the demodulation has a max amplification of 1), so stability is determined by the reduced model.

Let the map $L:\hat u \rightarrow \xi$ represent the static relation between $\hat u$ and $\xi$ in the reduced model. Assuming a constant $\hat u$, and that the amplitude $a$ is sufficiently small such that we may approximate the plant as locally linear with $G_{\hat u}(s)$, it is straight forward to show that $\xi$ is given by
\begin{align}
\xi = &\frac{a}{2}|F_H(i\omega)||G_{\hat u}(i\omega)||F_L(0)|\cos(\varphi_{\hat u}) + \\
&\frac{a}{2}|F_H(i\omega)||G_{\hat u}(i\omega)||F_L(i2\omega)|\cos (2\omega t +\varphi_{\hat u}).
\end{align}
where the phase-lag $\varphi_{\hat u} =\angle G_{\hat u}(i\omega) + \angle F_H(i\omega)$. 
Assume that the low-pass filter break-off frequency is low, i.e., $\omega_l \ll 1$ (but still larger than $k$ such that the timescale arguments are valid). We then get $|F_L(i2\omega)|\ll 1$ such that we can neglect\footnote{Alternatively we may apply averaging for the same purpose.} the high-frequency term to get
\begin{equation}
L(\hat u) = \frac{a}{2}|F_H(i\omega)||G_{\hat u}(i\omega)||F_L(0)|\cos(\varphi_{\hat u}).
\end{equation}
The reduced model is now
\begin{equation}
\dot {\hat u } = k L(\hat u).
\end{equation}
Note that the assumption that we operate at a stationary solution implies that we have $\varphi_{\hat u} = \pi/2 +n\pi$ for some $n\in\mathbf Z$, and consequently that $L(\hat u) = 0$, i.e., the stationarity condition \eqref{eq:stat_cond_lp} also holds for the reduced system.
We may determine stability of the above system by Lyapunov's indirect method, i.e., by linearizing $L$ at $\hat u$. The reduced model is thus stable when the sign of $k\cdot dL/du$ is negative. Using $\cos(\varphi_{\hat u}) = 0$ and $\sin(\varphi_{\hat u}) = -1^n$, we get
\begin{align}
\frac{dL}{d\hat u}=-1^n\frac{\varphi_{\hat u}}{d\hat u}\frac{a}{2}|F_H(i\omega)||G_{\hat u}(i\omega)||F_L(0)|.
\end{align}
Essentially any periodic solution which satisfy the necessary condition
\eqref{eq:stat_cond_lp} may hence be stabilized by the choice of $k$.

\section{Example: Classical ESC applied to a tubular reactor}
Consider a tubular isothermal reactor with plug flow used to convert $A$~$\rightarrow$~$B$, but where there is a side-reaction $B$~$\rightarrow$~$C$ incurring some loss of product. Due to the side-reaction, the reactor shows an optimum with respect to the residence time in the reactor which may be controlled via the feed-rate. Assuming standard mass-action kinetics, the dynamics of the reactor may be modeled by a system of PDEs of the form
\begin{equation}
\begin{aligned}
\frac{\partial a}{\partial t} &= -v\frac{\partial a}{\partial z} - k_1a, & a(t,0)&=a_0(t)\\
\frac{\partial b}{\partial t} &= -v\frac{\partial b}{\partial z} +k_1a - k_2b, & b(t,0)&=b_0(t)\\
\end{aligned}
\label{eq:reactor}
\end{equation}

\begin{figure}[htb!]
\includegraphics[]{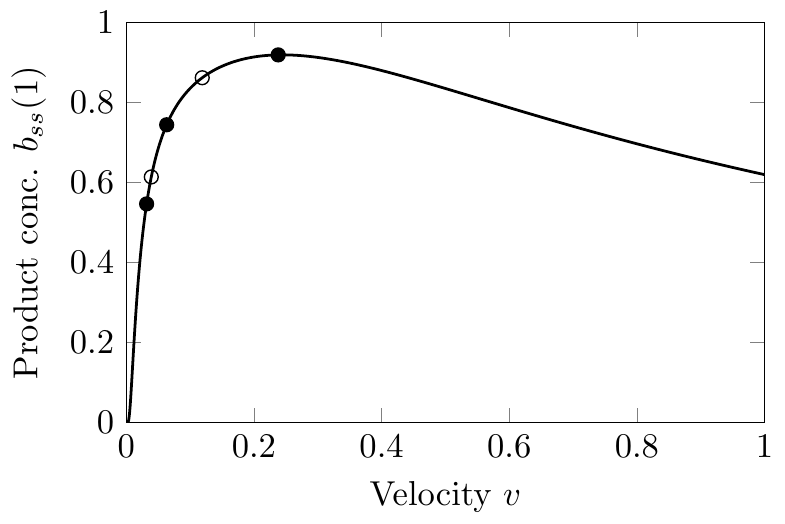}
\caption{Steady-state input-output map for the system \eqref{eq:reactor} with parameters $a_0 = 1, b_0 = 0, k_1 = 1, k_2 = 0.02$. The circles correspond to points where the stationarity condition \eqref{eq:stat_cond_lp} is satisfied using an ESC controller with parameters $\omega=0.4$, $\omega_l=\omega_h=0.1\omega$, $k=0.01$, $a=0.001$. The filled circles and the empty circles represent stable and unstable stationary solutions, respectively.}
\label{fig:inout}
\end{figure}

where $a$ and $b$ are dimensionless concentrations of $A$ and $B$ respectively, $k_i$ are the affinity constants, $z\in[0\, 1]$ is a normalized dimensionless space coordinate, and $v$ is the normalized velocity inside the reactor. Here $v$ is proportional to the feed-rate, and for simplicity we make use of $v$ directly as the control input. 
The above system of PDEs is converted to a system of ODEs using the method of lines by applying backward Euler over the spatial domain with $n=40$ discretization points.
Assume constant concentrations at the inlet and let the parameters be
\begin{equation}
a_0 = 1, \quad b_0 = 0,\quad k_1 =1, \quad k_2 = 0.02.
\end{equation}
This results in the steady-state input-output map depicted by the solid line in Figure~\ref{fig:inout}. 

Assume now that these dynamics are unknown and that classic ESC is applied in an attempt to optimize the production of $B$ with respect to the feed rate (i.e., the velocity $v$), using the ESC parameters
\begin{equation}
\omega = 0.4, \quad \omega_l=\omega_h = 0.1\omega, \quad k = 0.01, \quad a=0.001.
\end{equation}
Solving for condition \eqref{eq:stat_cond_lp} then yields the stationary solutions marked in Figure~\ref{fig:inout}.  The three solutions marked with filled circles are stable, and the solutions marked by empty circles are unstable.

\begin{figure}[htb!]
\includegraphics[]{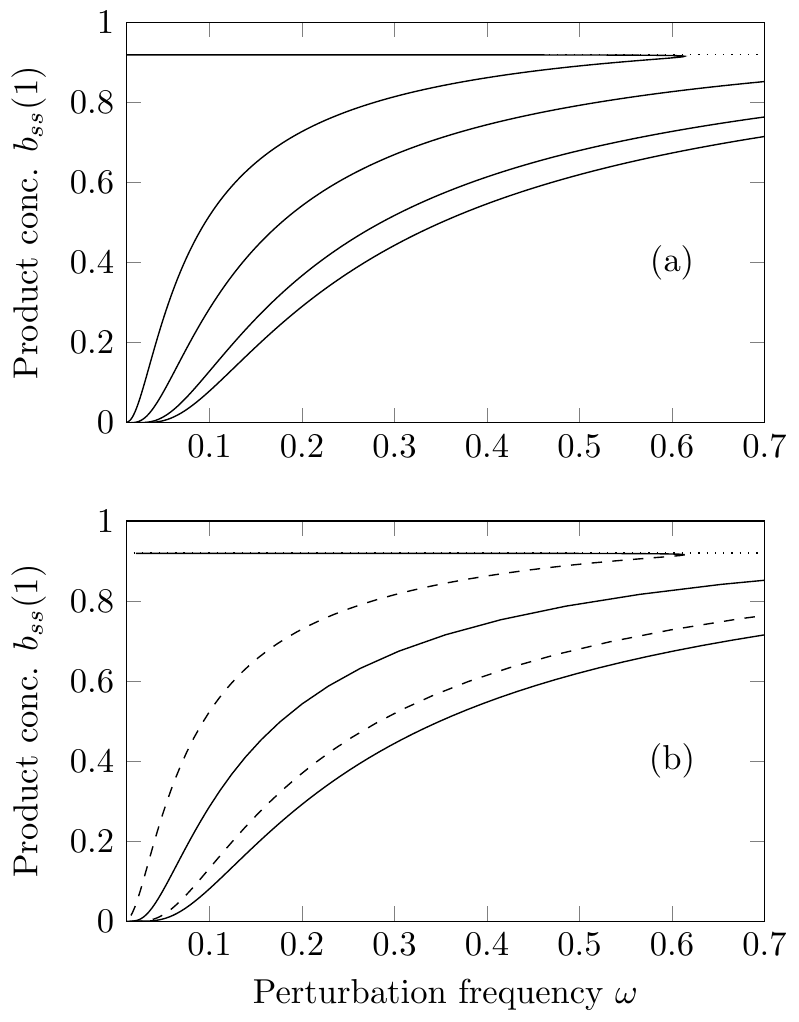}
\caption{Bifurcation diagrams for the closed-loop reactor with parameters $a_0 = $1, $b_0 = 0$, $k_1 = 1$, $k_2 = 0.02$, and ESC parameters $\omega_l=\omega_h=0.1\omega$, $k=0.01$, $a=0.001$, where $\omega$ is the varied parameter. The stationary solutions are represented by the steady-state output $b_{ss}$ at the end of the reactor. The diagram in (a) corresponds to stationary period-one solutions as determined by condition \eqref{eq:stat_cond_lp}. Stability is not considered in (a). The diagram in (b) corresponds to period-one solutions determined using AUTO \cite{Doedel2007}. In (b) the stability of the solutions is indicated by solid lines for stable limit cycles and dased lines for unstable limit cycles. The dotted line corresponds to the optimal output in both diagrams.  }
\label{fig:bifdiag}
\end{figure}

By solving equation \eqref{eq:stat_cond_lp} over a range of frequencies, a bifurcation diagram for the period-one stationary solutions to the ESC-loop may be determined, see Figure~\ref{fig:bifdiag}a. The bifurcation software AUTO\cite{Doedel2007} is used to verify this result; Figure~\ref{fig:bifdiag}b illustrates the period-one solutions determined using AUTO, solid for stable limit cycles and dashed for unstable limit cycles. 
As is clear from the figure, the near-optimal branch undergoes a fold bifurcation at the frequency $\omega=0.614$. This frequency cannot be considered low as compared to the time-constants of the local dynamics at the optimum, and this shows that the optimum may well be located and tracked using relatively fast estimation and control in the ESC-loop. However, the bifurcation also introduces local solution multiplicity and the second branch emerging from the fold bifurcation remains also for very low frequencies. As this solution branch is unstable, these solutions will not be observed in practice. However, the solution branch is still of practical interest as it acts as a separatrix, and hence limits the domain of attraction of the near-optimal solution branch. 
Furthermore, the unstable branch undergoes another fold bifurcation at a frequency close to zero, and this results in another stable branch with low product concentrations. As can be seen from Figure \ref{fig:bifdiag}, there are also two other solution branches for low product concentrations emerging from another fold bifurcation at low frequency. Thus, there are a total of 5 possible stationary solutions, of which 3 are stable, over a large frequency range. 
The fact that the solution branches remain as the frequency approaches zero illustrates that solution multiplicity may well occur also in cases where the loop is `conservatively' tuned. As a remark, we note that in e.g., \cite{krstic2000,tan2006}, it is assumed that the ESC-parameters may be selected such that we achieve time-scale separation between the plant and the ESC-controller. However, when the control input affects the time-constants of the plant, e.g., via the retention-time as in this example, such an assumption cannot be guaranteed to hold globally as the local dynamics of the plant may become arbitrarily slow in certain operating regions. Finally we remark that also other types of bifurcations than fold bifurcations may occur. For example, as $k$ is increased, some of the stable period-one solutions will loose stability in a period-doubling bifurcation, resulting in existence of a stable period-two solution.
 
\section{ Conclusions and discussion}
In Krstić and Wang \cite{krstic2000}, it is shown that the classical perturbation-based ESC method will possess a stable stationary solution in a local neighborhood of the optimum also when applied to systems with dynamics, provided the loop is conservatively tuned. However, existence of such a solution does not guarantee convergence in general; also uniqueness and the domain of attraction must be considered. In this paper, we approach the problem of uniqueness by considering conditions for stationarity without regard to the optimality conditions, and without resorting to asymptotic methods such as singular perturbations or averaging. We show that the local phase-lag is central for existence of stationary solutions by deriving a necessary condition which essentially any periodic stationary solution of the ESC-loop must satisfy.  Since it has been shown previously that a near-optimal solution exists \cite{krstic2000}, this result implies that there exists a connection between the local phase-lag and optimality. This connection is explored and it is shown that an extremum point in the steady-state input-output map generally corresponds to a bifurcation of the zero dynamics which in turn is reflected in large variations in the phase-lag locally about the optimum. These variations then ensure that the phase-lag condition will be satisfied locally about the optimum. However, the phase-lag condition may also be satisfied at points with no connection to the optimum whatsoever hence indicating that multiple stationary solutions can coexist. By applying elements of bifurcation theory to the necessary condition derived in the paper, we show that non-monotone variations in the local phase-lag with respect to the operating point is related to existence of fold bifurcations, a common source of solution multiplicity. Given the existence of a stationary solution, a simplified stability analysis is provided to show that essentially any of the stationary solutions discussed in the paper in principle may be stabilized by the choice of the integral gain $k$. Finally, a simple example is provided to illustrate the results.

It has long been recognized that the phase-lag affects the performance of classic ESC \cite{sternby1979}. However, the effect of the phase has mainly been considered in relation to stability and the convergence rate of the scheme. Here we show that the phase in fact is instrumental for successful operation of the loop in the dynamic case, not least in that it is critical for the existence of stationary solutions, including the near-optimal solution. 

That the classical ESC-method may display solution multiplicity may have severe implications for the applicability of the method as a general purpose optimization technique; if the near-optimal solution is not unique, convergence to this solution cannot be guaranteed other than from operating points in a local neighborhood of the optimum. Note that this holds regardless of the stability properties of the solution(s) that are unrelated to the optimum, that is, also the existence of additional unstable solutions may limit the domain of attraction of the near-optimal solution.

In the example presented in this paper, we first observe that the near-optimal branch remains close to the optimum also for relatively aggressive tuning of the loop. Even though the domain of attraction in this case becomes relatively small, such a tuning provides a higher bandwidth and may hence be of use for tracking an optimum that moves over time, even though it may be less useful for locating the optimum initially since this would require a very accurate initial guess for convergence. Second, we observe that the near-optimal solution branch ends in a fold bifurcation from which also an unstable solution branch, unrelated to the optimum, emerges. This branch continues to exist also for relatively low frequencies, hence implying that solution multiplicity remains also in parameter domains where the loop is conservatively tuned. 

\bibliographystyle{plain}
\bibliography{references}

\end{document}